\documentstyle[12pt]{article}
\input amssym.def
\begin{document}
\def \Z{\Bbb Z}
\def \C{\Bbb C}
\def \R{\Bbb R}
\def \Q{\Bbb Q}
\def \N{\Bbb N}
\def \P{\Bbb P}
\def \bZ{{\bf Z}}
\def \TZ{{1\over T}{\bf Z}}
\def \wt{{\rm wt}\;}
\def \g{\frak g}
\def \mod{{\rm mod}\;}
\def \pf{{\bf Proof}.}
\def \x{{\bf x}}
\def \y{{\bf y}}

\def \wt{{\rm wt}}
\def \span{{\rm span}}
\def \Res{{\rm Res}}
\def \End{{\rm End}}
\def \Hom{{\rm Hom}}
\def \<{\langle} 
\def \>{\rangle}
\def \be{\begin{equation}\label}
\def \ee{\end{equation}}
\def \bex{\begin{exa}\label}
\def \eex{\end{exa}}
\def \bl{\begin{lem}\label}
\def \el{\end{lem}}
\def \bt{\begin{thm}\label}
\def \et{\end{thm}}
\def \bp{\begin{prop}\label}
\def \ep{\end{prop}}
\def \br{\begin{rem}\label}
\def \er{\end{rem}}
\def \bc{\begin{coro}\label}
\def \ec{\end{coro}}
\def \bd{\begin{de}\label}
\def \ed{\end{de}}

\newtheorem{thm}{Theorem}[section]
\newtheorem{prop}[thm]{Proposition}
\newtheorem{coro}[thm]{Corollary}
\newtheorem{conj}[thm]{Conjecture}
\newtheorem{exa}[thm]{Example}
\newtheorem{lem}[thm]{Lemma}
\newtheorem{rem}[thm]{Remark}
\newtheorem{de}[thm]{Definition}
\newtheorem{hy}[thm]{Hypothesis}
\makeatletter
\@addtoreset{equation}{section}
\def\theequation{\thesection.\arabic{equation}}
\makeatother
\makeatletter

\begin{center}
{\Large \bf Simple vertex operator algebras are nondegenerate}
\end{center}
\begin{center}{Haisheng Li\footnote{Partially supported by NSF grant
DMS-9970496 and a grant from Rutgers Research Council}\\
Department of Mathematical Sciences, Rutgers University, Camden, NJ 08102\\
and\\
Department of Mathematics, Harbin Normal University, Harbin, China}
\end{center}

\begin{abstract}
In this paper it is shown that every irreducible vertex algebra of
countable dimension and every simple vertex operator algebra
is nondegenerate in the sense of Etingof and Kazhdan.
\end{abstract}

\section{Introduction}

In one of their series papers on quantization of Lie bialgebras, 
Etingof and Kazhdan [EK] introduced a notion of
nondegeneracy of a vertex operator algebra.
A vertex operator algebra $V$ is said to be nondegenerate
if for every positive integer $n$ the linear map $Z_{n}$ from
$V^{\otimes n}\otimes \C((z_{1},\dots,z_{n}))$ to $V((z_{1}))\cdots
((z_{n}))$ defined by
$$Z_{n}(v^{(1)}\otimes \cdots \otimes v^{(n)}\otimes f)=
fY(v^{(1)},z_{1})\cdots Y(v^{(n)},z_{n}){\bf 1}$$
is injective.
It was proved therein that in their definition of the notion 
of braided vertex operator algebra, if the classical limit 
vertex operator algebra is nondegenerate, certain
axioms are automatically satisfied.
Let $V(\g,K)$ be the vertex operator algebra associated with 
a (finite-dimensional) Lie algebra $\g$ and a (complex) number $K$
(cf. [FF], [FZ], [Lia]).
It was proved in [EK] that if $V(\g,K)$ is an irreducible 
$\hat{\g}$-module, then $V(\g,K)$ is nondegenerate. Notice that
the irreducibility of $V(\g,K)$ as a $\hat{\g}$-module amounts to
the simplicity of the vertex operator algebra  $V(\g,K)$.
So, it is reasonable to conjecture that general simple vertex operator
algebras are nondegenerate. In this paper, we shall prove that
this conjecture is indeed true.

In literature, there are certain results closely related to
the injectivity of the maps $Z_{n}$. In [DL], among other results
it was proved that if $V$ is a simple vertex
operator algebra and $W$ is an irreducible $V$-module, then
$Y(v,z)w\ne 0$ for $0\ne v\in V,\; 0\ne w\in W$.
Furthermore, among other results it was proved in [DM] that 
the linear map $Y$ viewed as a map from $V\otimes W$ to $W((z))$ 
is injective. 
It is quite conceivable that ideas in [DL], [DM] and [EK]
will be highly valuable for proving the conjecture. 
Indeed, our proof here essentially uses the same ideas.

Notice that in the notion of vertex operator algebra used in [EK],
no conformal vector and no $\Z$-grading are assumed. A vertex operator 
algebra in the sense of [EK] is often called a vertex algebra.
In this paper we consider a general vertex algebra $V$.
For any $V$-module $W$ we define a linear map $Z_{n}^{W}$ from 
$V^{\otimes n}\otimes W\otimes \C((z_{1}))\cdots ((z_{n}))$ to
$W((z_{1}))\cdots ((z_{n}))$ by
$$Z_{n}^{W}(v^{(1)}\otimes \cdots \otimes v^{(n)}\otimes w\otimes f)=
fY(v^{(1)},z_{1})\cdots Y(v^{(n)},z_{n})w.$$
It follows from [FHL] that
$V^{\otimes (n+1)}\otimes \C((z_{1}))\cdots ((z_{n}))$
is a natural vertex algebra with 
$V^{\otimes n}\otimes W\otimes \C((z_{1}))\cdots ((z_{n}))$
as a module.
We first prove that $\ker Z_{n}^{W}$ is a submodule.
To describe submodules of 
$V^{\otimes n}\otimes W\otimes \C((z_{1}))\cdots ((z_{n}))$
we slightly generalize the result of [FHL] on the irreducibility of
tensor product modules for tensor product vertex operator
algebras in the context of vertex algebras of countable dimension.
Using this, we show that if $V$ has countable dimension and irreducible
in the sense that $V$ is an irreducible $V$-module
and if $W$ is a $V$-module, then
any submodule of 
$V^{\otimes n}\otimes W\otimes \C((z_{1}))\cdots ((z_{n}))$
is of the form $V^{\otimes n}\otimes U$, where $U$ is 
a $V\otimes \C((z_{1}))\cdots ((z_{n}))$-submodule of
$W\otimes \C((z_{1}))\cdots ((z_{n}))$. Then it follows that
$\ker Z_{n}^{W}=0$. From this we show that
every irreducible vertex algebra of countable dimension
is nondegenerate. In particular, this implies that
every simple vertex operator algebra in the sense of [FLM] and [FHL]
is nondegenerate.

\section{The main result}
We here shall recall the notion of nondegeneracy of a vertex 
algebra from [EK] and prove that every irreducible vertex algebra
of countable dimension is nondegenerate.

Throughout this paper, the ground field is $\C$.
In this paper we shall use the following definition of the notion of 
vertex algebra (cf. [B], [FLM], [FHL], [DL], [Li2], [LL]):

\bd{dvadefinition}
{\em A {\em vertex algebra} is a vector space $V$ equipped with 
a linear map, called the {\em vertex operator map,}
\begin{eqnarray}
Y: V &\rightarrow& (\End V)[[z,z^{-1}]]\nonumber\\
v &\mapsto & Y(v,z)=\sum_{n\in \Z}v_{n}z^{-n-1}\;\;\;( v_{n}\in \End V)
\end{eqnarray}
and equipped with a distinguished vector ${\bf 1}\in V$,
called the {\em vacuum vector},
such that the following axioms hold: For $u,v\in V$, 
\begin{eqnarray}\label{etruncationva}
u_{n}v=0\;\;\;\mbox{ for $n$ sufficiently large};
\end{eqnarray}
\begin{eqnarray}\label{evacuumva}
Y({\bf 1},z)=1;
\end{eqnarray}
for $v\in V$, 
\begin{eqnarray}
Y(v,z){\bf 1}\in V[[z]]\;\;\;\mbox{ and }\;\;\; 
Y(v,z){\bf 1}|_{z=0}\left(=v_{-1}{\bf 1}\right)=v;
\end{eqnarray}
and for $u,v\in V$,
\begin{eqnarray}\label{ejacobiva}
& &z_{0}^{-1}\delta\left(\frac{z_{1}-z_{2}}{z_{0}}\right)
Y(u,z_{1})Y(v,z_{2})-z_{0}^{-1}\delta\left(\frac{z_{2}-z_{1}}{-z_{0}}\right)
Y(v,x_{2})Y(u,x_{1})\nonumber\\
&=&z_{2}^{-1}\delta\left(\frac{z_{1}-z_{0}}{z_{2}}\right)
Y(Y(u,z_{0})v,z_{2})
\end{eqnarray}
(the {\em Jacobi identity}).}
\ed

\br{requivalence}
{\em In literature, there are variant definitions
of the notion of vertex algebra (cf. [B], [DL], [Li2], [Ka]). 
It was proved in [Li2] that the definition 
given in [B] is equivalent to the current definition (with ground
field $\C$). It follows immediately from 
a known duality result ([G], [FHL], [DL], [Li2]) that 
the definition given in [EK] and [Ka] 
(where the vacuum vector is denoted by $\Omega$) 
is also equivalent to the current definition.}
\er

For a vertex algebra $V$, the vertex operator map
$Y$ is a linear map from $V$ to $\Hom (V,V((z)))$.
The map $Y$ can alternatively be considered as a 
linear map from $V\otimes V$ to $V((z))$.
Following [EK], we alternatively denote this map by $Y(z)$.
Define a linear operator ${\cal{D}}\in \End V$ by
\begin{eqnarray}
D(v)=v_{-2}{\bf 1}\left(=\left({d\over dz}Y(v,z){\bf 1}\right)|_{z=0}\right).
\end{eqnarray}
Then (cf. [LL])
\begin{eqnarray}
& &[{\cal{D}},Y(u,z)]=Y({\cal{D}}u,z)={d\over dz}Y(u,z),\\
& &Y(u,z)v=e^{z{\cal{D}}}Y(v,-z)u\;\;\;\mbox{ for }u,v\in V.
\end{eqnarray}

It was proved ([FLM], [FHL], [DL], [Li2]) that the Jacobi identity
is equivalent to the following {\em weak commutativity and associativity}:
For any $u,v\in V$, there exists a nonnegative integer $k$ such that
\begin{eqnarray}\label{eweakcomm}
(z_{1}-z_{2})^{k}Y(u,z_{1})Y(v,z_{2})=(z_{1}-z_{2})^{k}Y(v,z_{2})Y(u,z_{1});
\end{eqnarray}
and for any $u,v,w\in V$ there exists a nonnegative integer $l$ such that
\begin{eqnarray}\label{eweakasso}
(z_{0}+z_{2})^{l}Y(u,z_{0}+z_{2})Y(v,z_{2})w
=(z_{0}+z_{2})^{l}Y(Y(u,z_{0})v,z_{2})w.
\end{eqnarray}

A {\em $V$-module} (cf. [Li2]) is a vector space $W$ equipped with
a linear map $Y$ from $V$ to $(\End W)[[z,z^{-1}]]$ such that all the axioms
in defining the notion of vertex algebra that makes sense hold. That is,
the truncation condition (\ref{etruncationva}), 
the vacuum property (\ref{evacuumva})
and the Jacobi identity (\ref{ejacobiva}) hold.

The notion of ideal is defined in the obvious way; an ideal of a
vertex algebra $V$ is
a subspace $U$ such that $u_{m}v,\;v_{m}u\in U$ for all $v\in V,\; u\in U$.

\bd{dsimpleva}
{\em A vertex algebra $V$ is said to be {\em simple} if there is no
nontrivial ideal and $V$ is said to be {\em irreducible}
if $V$ is an irreducible $V$-module.}
\ed

Clearly, an irreducible vertex
algebra is simple, but simple vertex algebras are
not necessarily irreducible. For example, the vertex
algebra constructed in [B] from the commutative associative algebra
$\C[x]$ with the standard derivation is simple, but not irreducible.
On the other hand, simple vertex operator algebras 
in the sense of [FLM] and [FHL] are always
irreducible because any submodule of $V$ is an ideal (cf. [FHL]).

{\em Let $V$ be a vertex algebra, fixed throughout this section.}
Following [EK], for a positive integer $n$,
we define a linear map $Z_{n}$ from 
$V^{\otimes n}\otimes \C((z_{1},\dots,z_{n}))$
to $V((z_{1}))\cdots ((z_{n}))$ by
\begin{eqnarray}
Z_{n}(v^{(1)}\otimes \cdots\otimes v^{(n)}\otimes f)
=f Y(v^{(1)},z_{1})\cdots Y(v^{(n)},z_{n}){\bf 1}.
\end{eqnarray}
In [EK], $Z_{n}$ was alternatively defined as
\begin{eqnarray}\label{edefinezn}
Z_{n}=Y(z_{1})(1\otimes Y(z_{2}))\cdots (1^{\otimes n-1}\otimes Y(z_{n}))
(1^{\otimes n}\otimes {\bf 1}).
\end{eqnarray}
The following notion was due to [EK]:

\bd{dnondegeneracy}
{\em A vertex algebra $V$ is said to be {\em nondegenerate} if
the linear maps $Z_{n}$ are injective for all positive integers $n$.}
\ed

\br{rcase1}
{\em Consider the case $n=1$. For $v\in V,\; f\in \C((z))$, we have
$$Z_{1}(v\otimes f)=fY(v,z){\bf 1}=fe^{z{\cal{D}}}Y({\bf 1},z)v
=e^{z{\cal{D}}}fv.$$
Then $Z_{1}=e^{z{\cal{D}}}\pi$, where $\pi$ is the natural embedding of
$V\otimes \C((z))$ into $V((z))$.
It follows immediately that $Z_{1}$ is injective.}
\er

Now, let $W$ be a $V$-module. Then
we similarly define a linear map $Z_{n}^{W}$ from
$V^{\otimes n}\otimes W\otimes \C((z_{1}))\dots ((z_{n}))$ to 
$W((z_{1}))\cdots ((z_{n}))$ by
\begin{eqnarray}
Z_{n}^{W}(v^{(1)}\otimes \cdots \otimes v^{(n)}\otimes w\otimes f)
=f Y(v^{(1)},z_{1})\cdots Y(v^{(n)},z_{n})w
\end{eqnarray}
for $v^{(1)},\dots,v^{(n)}\in V,\; w\in W,\; f\in \C((z_{1}))\dots ((z_{n}))$.
Notice that 
$\C((z_{1}))\cdots ((z_{n}))$ is a commutative associative algebra (with
identity element) and that for any vector space $U$,
$U((z_{1}))\cdots ((z_{n}))$ is a $\C((z_{1}))\cdots ((z_{n}))$-module.
It is clear that
$Z_{n}^{W}$ is $\C((z_{1}))\dots ((z_{n}))$-linear where
$V^{\otimes n}\otimes W\otimes \C((z_{1}))\dots ((z_{n}))$
is considered as a $\C((z_{1}))\dots ((z_{n}))$-module in the obvious way.
Let $E_{n}$ be the embedding of $V^{\otimes n}$ into $V^{\otimes (n+1)}$
defined by
\begin{eqnarray}
E_{n}(v^{(1)}\otimes\cdots\otimes v^{(n)})=
v^{(1)}\otimes\cdots\otimes v^{(n)}\otimes {\bf 1}.
\end{eqnarray}
Then
\begin{eqnarray}
Z_{n}=Z_{n}^{V}E_{n}.
\end{eqnarray}

In literature, there are certain results which are closely related
to the injectivity of linear maps $Z_{n}^{V}$.
The following result is due to Dong and Mason [DM] while
the special case is due to Dong and Lepowsky [DL]:

\bp{pdmdl}
Let $V$ be a simple vertex operator algebra and let $W$ be an
irreducible $V$-module. Let $v^{(1)},\dots,
v^{(r)}\in V$ be nonzero vectors and let
$w^{(1)},\dots,w^{(r)}\in W$ be linearly independent vectors. Then 
\begin{eqnarray}
\sum_{i=1}^{r}Y(v^{(i)},z)w^{(i)}\ne 0.
\end{eqnarray}
In particular,
\begin{eqnarray}
Y(v,z)w\ne 0\;\;\;\mbox{ for }0\ne v\in V,\; 0\ne w\in W.
\end{eqnarray}
\ep

\br{rdm}
{\em Proposition \ref{pdmdl} exactly asserts that 
$Y$ viewed as a linear map from $V\otimes W$ to $W((z))$ is injective.}
\er

{}From [FHL], for any positive integer $n$, $V^{\otimes (n+1)}$ has a
natural vertex algebra structure and for $V$-modules
$W_{1},\dots, W_{n+1}$, $W_{1}\otimes \cdots\otimes W_{n+1}$ 
has a natural $V^{\otimes (n+1)}$-module structure.
For $1\le i\le n+1$, let $\pi_{i}$ be the embedding of $V$ into
$V^{\otimes (n+1)}$ sending $v$ to the product vector all of whose 
factors but the $i$th factor which is $v$ are the vacuum vector ${\bf 1}$.
Let $U$ be any subspace of $W_{1}\otimes \cdots\otimes W_{n+1}$. 
It is clear (cf. [FHL]) that 
$U$ is a submodule if and only
if for $i=1,\dots,n+1$,
\begin{eqnarray}
Y(\pi_{i}(v),z)U\subset U((z))\;\;\;\mbox{ for all }v\in V.
\end{eqnarray}

Note that any commutative associative algebra with identity is 
a natural vertex algebra. Then from [FHL], 
$V^{\otimes (n+1)}\otimes \C((z_{1}))\cdots ((z_{n}))$ is a vertex algebra. 
Now we state our first result:

\bp{pkernal}
For any positive integer $n$, $\ker Z_{n}^{W}$
is a $V^{\otimes (n+1)}\otimes \C((z_{1}))\dots ((z_{n}))$-submodule
of $V^{n}\otimes W\otimes \C((z_{1}))\dots ((z_{n}))$.
\ep

For convenience we first prove the following simple fact:

\bl{lsimplefact}
Let $U$ be a vector space and let
\begin{eqnarray}
f(z_{1},\dots,z_{n})\in U((z_{1}))\cdots ((z_{n})).
\end{eqnarray}
Assume that there exist nonnegative integers $k_{ij}$ 
for $1\le i<j\le n$ such that
\begin{eqnarray}\label{eassumption}
\left(\prod_{1\le i<j\le n}(z_{i}-z_{j})^{k_{ij}}\right)
f(z_{1},\dots,z_{n})=0.
\end{eqnarray}
Then $f(z_{1},\dots,z_{n})=0$.
\el

{\bf Proof.} The key issue here is the cancelation law. From [FLM],
for any three formal series $A, B$ and $C$, if $ABC$, $AB$ and $BC$
all exist (algebraically), then 
\begin{eqnarray}
A(BC)=(AB)C=ABC.
\end{eqnarray}
Set 
$$A=\prod_{1\le i<j\le n}(z_{i}-z_{j})^{-k_{ij}},\;\;\;
B=\prod_{1\le i<j\le n}(z_{i}-z_{j})^{k_{ij}},\;\;\;\; 
C=f(z_{1},\dots,z_{n}),$$
where we use the usual binomial expansion convention (cf. [FLM]). 
Then it follows immediately from (\ref{eassumption})
that $f(z_{1},\dots,z_{n})=0$.
$\;\;\;\;\Box$

{\bf Proof of Proposition \ref{pkernal}.} 
Since $Z_{n}^{W}$ is already $\C((z_{1}))\cdots ((z_{n}))$-linear,
what we must prove is that for $i=1,\dots,n+1$,
\begin{eqnarray}
Y(\pi_{i}(v),z)Z_{n}^{W}\subset Z_{n}^{W}((z))\;\;\;\mbox{ for }v\in V.
\end{eqnarray}
That is, we must prove that if 
$X\in V^{\otimes n}\otimes W\otimes \C((z_{1}))\cdots ((z_{n}))$ with
$Z_{n}^{W}(X)=0$, then for $i=1,\dots,n+1$,
\begin{eqnarray}\label{establemain}
Z_{n}^{W}(Y(\pi_{i}(v),z)X)=0\;\;\;\mbox{ for all }v\in V.
\end{eqnarray}
We shall prove this in four steps.

Claim 1: For $i=n+1$, (\ref{establemain}) holds.

Let 
\begin{eqnarray}\label{eelementx}
X=\sum_{\alpha_{i}, \beta, \gamma}
u^{1\alpha_{1}}\otimes \cdots \otimes u^{n\alpha_{n}}\otimes
w^{\beta}\otimes f_{\gamma},
\end{eqnarray}
(a {\em finite} sum). We shall use this general $X$ for the whole proof.
Then 
\begin{eqnarray}\label{eznx=0}
\sum_{\alpha_{i}, \beta, \gamma}
f_{\gamma} Y(u^{1\alpha_{1}},z_{1})\cdots Y(u^{n\alpha_{n}},z_{n})w^{\beta}
=Z_{n}^{W}(X)=0.
\end{eqnarray}
Let $v\in V$. From the weak commutativity there exists 
a nonnegative integer $k$ such that
\begin{eqnarray}\label{e2.9weakcomm}
(z-z_{i})^{k}Y(v,z)Y(u^{i\alpha_{i}},z_{i})
=(z-z_{i})^{k}Y(u^{i\alpha_{i}},z_{i})Y(v,z)
\end{eqnarray}
for all the indices $i$ and $\alpha_{i}$.
Multiplying (\ref{eznx=0}) 
by $(z-z_{1})^{k}\cdots (z-z_{n})^{k}Y(v,z)$ (from left)
and then using (\ref{e2.9weakcomm}) we get
\begin{eqnarray}
(z-z_{1})^{k}\cdots (z-z_{n})^{k}\sum_{\alpha_{i}, \beta, \gamma}
f_{\gamma} Y(u^{1\alpha_{1}},z_{1})\cdots Y(u^{n\alpha_{n}},z_{n})
Y(v,z)w^{\beta}=0.
\end{eqnarray}
In view of Lemma \ref{lsimplefact} we have
\begin{eqnarray}
\sum_{\alpha_{i}, \beta, \gamma}
f_{\gamma} Y(u^{1\alpha_{1}},z_{1})\cdots Y(u^{n\alpha_{n}},z_{n})
Y(v,z)w^{\beta}=0.
\end{eqnarray}
Thus $Z_{n}^{W}(Y(\pi _{n+1}(v),z)X)=0$. 

Claim 2: For $i=n$, (\ref{establemain}) holds. From (\ref{eznx=0}),
using (\ref{e2.9weakcomm}) we get
\begin{eqnarray}
\left(\prod_{i=1}^{n-1}(z-z_{i})^{k}\right)
\sum _{\alpha_{i}, \beta, \gamma}
f_{\gamma} Y(u^{1\alpha_{1}},z_{1})\cdots Y(u^{(n-1)\alpha_{n-1}},z_{n-1})
Y(v,z)Y(u^{n\alpha_{n}},z_{n})w^{\beta}=0.
\end{eqnarray}
In view of Lemma \ref{lsimplefact} we have
\begin{eqnarray}
\sum _{\alpha_{i}, \beta, \gamma}
f_{\gamma} Y(u^{1\alpha_{1}},z_{1})\cdots Y(u^{n-1\alpha_{n-1}},z_{n-1})
Y(v,z)Y(u^{n\alpha_{n}},z_{n})w^{\beta}=0.
\end{eqnarray}
{}From the weak associativity, there exists
a nonnegative integer $l$ such that
\begin{eqnarray}
(z_{0}+z_{n})^{l}Y(v,z_{0}+z_{n})Y(u^{n\alpha_{n}},z_{n})w^{\beta}
=(z_{0}+z_{n})^{l}Y(Y(v,z_{0})u^{n\alpha_{n}},z_{n})w^{\beta}
\end{eqnarray}
for all the indices $\alpha_{n}$ and $\beta$, (which are finitely many).
Then
\begin{eqnarray}
& &(z_{0}+z_{n})^{l}
\sum _{\alpha_{i}, \beta, \gamma}
f_{\gamma} Y(u^{1\alpha_{1}},z_{1})\cdots Y(u^{(n-1)\alpha_{n-1}},z_{n-1})
Y(Y(v,z_{0})u^{n\alpha_{n}},z_{n})w^{\beta}\nonumber\\
&=&(z_{0}+z_{n})^{l}
\sum _{\alpha_{i}, \beta, \gamma}
f_{\gamma} Y(u^{1\alpha_{1}},z_{1})\cdots Y(u^{(n-1)\alpha_{n-1}},z_{n-1})
Y(v,z_{0}+z_{n})Y(u^{n\alpha_{n}},z_{n})w^{\beta}\nonumber\\
&=&0.
\end{eqnarray}
By Lemma \ref{lsimplefact} we get
\begin{eqnarray}
\sum _{\alpha_{i}, \beta, \gamma}
f_{\gamma} Y(u^{1\alpha_{1}},z_{1})\cdots Y(u^{(n-1)\alpha_{n-1}},z_{n-1})
Y(Y(v,z_{0})u^{n\alpha_{n}},z_{n})w^{\beta}=0.
\end{eqnarray}
That is, $Z_{n}^{W}(Y(\pi_{n}(v),z)X)=0$.

Claim 3: $\ker Z_{n}^{W}$ is stable under the action of 
$S_{n}$. 
(Note that the symmetric group $S_{n}$ naturally acts 
on $V^{\otimes n}\otimes W\otimes \C((z_{1}))\cdots ((z_{n}))$.)

It suffices to prove that $\ker Z_{n}^{W}$ is stable under
the actions of transitions $\sigma_{i(i+1)}$ for $i=1,\dots, n-1$.
Let $X\in \ker Z_{n}^{W}$ and write $X$ as in (\ref{eelementx}).
{}From the weak commutativity there exists a
nonnegative integer $k$ such that
\begin{eqnarray}
(z_{i}-z_{i+1})^{k}Y(u^{i\alpha_{i}},z_{i})Y(u,z_{i+1})
=(z_{i}-z_{i+1})^{k}Y(u,z_{i+1})Y(u,z_{i})
\end{eqnarray}
for all the indices $i$ and $\alpha_{i}$. Then
\begin{eqnarray}
(z_{i}-z_{i+1})^{k}Z_{n}^{W}(\sigma_{i(i+1)}X)=
(z_{i}-z_{i+1})^{k}Z_{n}^{W}(X)=0.
\end{eqnarray}
In view of Lemma \ref{lsimplefact},
we have $Z_{n}^{W}(\sigma_{i(i+1)}X)=0$. That is, $\ker Z_{n}^{W}$ is
stable under the action of $\sigma_{i(i+1)}$ for $i=1,\dots,n$.
Therefore $\ker Z_{n}^{W}$ is $S_{n}$-stable.

Claim 4: For $1\le i\le n$, (\ref{establemain}) holds.

Let $\sigma$ be the permutation $(in)$. Then we have
\begin{eqnarray}
Y(\pi_{i}(v),z)X=\sigma(Y(\pi_{n}(v),z)\sigma(X))
\;\;\;\mbox{ for }v\in V.
\end{eqnarray}
By Claim 3, $\sigma(X)\in \ker Z_{n}^{W}$.
Furthermore, by Claim 2 we have
$$Y(\pi_{n}(v),z)\sigma(X)\in (\ker Z_{n}^{W})((z)).$$
By Claim 3 again we get
$$\sigma(Y(\pi_{n}(v),z)\sigma(X))\in (\ker Z_{n}^{W})((z)).$$
Therefore 
$$Z_{n}^{W}(Y(\pi_{i}(v),z)X)=Z_{n}^{W}(\sigma(Y(\pi_{n}(v),z)\sigma(X)))
=0.$$
This completes the proof.
$\;\;\;\;\Box$

Noticing that for $0\ne w\in W,\; 0\ne f\in \C((z_{1}))\dots ((z_{n}))$,
\begin{eqnarray}
Z_{n}^{W}({\bf 1}^{\otimes n}\otimes w\otimes f)=
f Y({\bf 1},z_{1})\cdots Y({\bf 1},z_{n})w=f w\ne 0,
\end{eqnarray}
we have
\begin{eqnarray}
\ker Z_{n}^{W}\ne V^{\otimes n}\otimes W\otimes \C((z_{1}))\dots ((z_{n})).
\end{eqnarray}
If $V^{n}\otimes W\otimes \C((z_{1}))\dots ((z_{n}))$ is
an irreducible $V^{\otimes (n+1)}\otimes
\C((z_{1}))\dots ((z_{n}))$-module, in view of Proposition \ref{pkernal}
we will immediately have $\ker Z_{n}^{W}=0$.
When $V$ is a simple vertex operator algebra and 
$W$ is an irreducible $V$-module, it follows from
[FHL] that $V^{\otimes n}\otimes W$ is an irreducible
$V^{\otimes (n+1)}$-module. However, the result of [FHL] does not
imply that 
$V^{\otimes n}\otimes W\otimes \C((z_{1}))\dots ((z_{n}))$
is an irreducible 
$V^{\otimes (n+1)}\otimes \C((z_{1}))\dots ((z_{n}))$-module.
Motivated by this, we prove the following simple result
(which might have been known somewhere):

\bl{ldensity}
Let $A_{1}$ and $A_{2}$ be associative algebras with identity element
and let $U_{1}$ and 
$U_{2}$ be modules for $A_{1}$ and $A_{2}$, respectively.
Assume that $U_{1}$ is irreducible and $\End _{A_{1}}U_{1}=\C$.
Then any $A_{1}\otimes A_{2}$-submodule of $U_{1}\otimes U_{2}$
is of the form $U_{1}\otimes U_{2}'$, where $U_{2}'$ is an
$A_{2}$-submodule of $U_{2}$. Furthermore, if $U_{2}$ is irreducible, 
then $U_{1}\otimes U_{2}$ is irreducible.
\el

{\bf Proof.} Fix a basis $\{ u_{1\alpha}\;|\; \alpha\in I\}$ for $U_{1}$.
Let $U$ be a nonzero $A_{1}\otimes A_{2}$-submodule of
$U_{1}\otimes U_{2}$. Let $u$ be any nonzero element of $U$. Then
$$u=u_{1\alpha_{1}}\otimes u_{21}+\cdots +u_{1\alpha_{r}}\otimes u_{2r},$$
where $u_{21},\dots,u_{2r}$ are (finitely many) nonzero vectors in $U_{2}$.
Since $U_{1}$ is irreducible and $\End _{A_{1}}U_{1}=\C$, 
in view of the density theorem
(cf. [J]), there exists $a\in A_{1}$ such that
$$au_{1\alpha_{1}}\ne 0\;\;\;\mbox{ and }\;\; au_{1\alpha_{i}}=0
\;\;\;\mbox{ for }i=2,\dots,r.$$
Then
$$0\ne au_{1\alpha_{1}}\otimes u_{21}=(a\otimes 1)u\in U.$$
Since $U_{1}$ is irreducible, we have $A_{1}au_{1\alpha_{1}}=U_{1}$.
Thus $U_{1}\otimes u_{21}=(A_{1}a\otimes 1)u\subset U$.
Similarly, we have $U_{1}\otimes u_{2i}\subset U$ for $i=2,\dots,r$.
Set
$$U_{2}'(u)=A_{2}u_{21}+\cdots +A_{2}u_{2r}\subset U_{2}.$$
(Notice that $u_{2i}$'s are uniquely determined by $u$.)
Then 
$$u\in U_{1}\otimes U_{2}'(u)\subset U.$$
Setting
\begin{eqnarray}
U_{2}'=\sum_{0\ne u\in U}U_{2}'(u)\subset U_{2},
\end{eqnarray}
we get $U=U_{1}\otimes U_{2}'$, completing the proof.
$\;\;\;\;\Box$

Closely related to Lemma \ref{ldensity} is
the following result which can be found in [Di], or [CG]:

\bl{lschur}
Let $A$ be an associative algebra with identity and 
let $U$ be an irreducible $A$-module of countable dimension.
Then $\End _{A}U=\C$.
\el


We now apply Lemma \ref{lschur} in vertex algebra theory.

\bl{lvacounatble}
Let $V$ be a vertex algebra of countable dimension and let $W$ be an
irreducible $V$-module. Then $W$ is of countable dimension and 
$\End _{V}W=\C$. 
\el

{\bf Proof.} Let $w$ be any nonzero element of $W$.
It was proved in  [DM] and [Li1] that the linear span of 
vectors $v_{m}w$ for $v\in V,\; m\in \Z$ is a submodule of $W$.
Consequently,
$$W={\rm span} \{ v_{m}w\;|\;v\in V,\; m\in \Z\}.$$
Then we immediately see that $W$ has countable dimension.
Furthermore, let $A$ be the subalgebra of $\End W$ generated
by all the operators $v_{m}$ for $v\in V,\; m\in \Z$. Then $A$ acts
irreducibly on $W$ and $\End_{A}W=\End_{V}W$.
In view of Lemma \ref{lschur} we have
$\End _{A}W=\C$. Thus $\End_{V}W=\End _{A}W=\C$.
$\;\;\;\;\Box$

The following result slightly generalizes the corresponding result of [FHL]:

\bp{pgeneralizationfhl}
Let $V_{1},\dots, V_{r}$ be vertex algebras of countable dimension 
and let $W_{1}$,$\dots$, $W_{r}$ be irreducible modules for
$V_{1},\dots,V_{r}$, respectively. Then
$W_{1}\otimes \cdots \otimes W_{r}$ is an irreducible $V_{1}\otimes
\cdots \otimes V_{r}$-module with
$$\End_{V_{1}\otimes\cdots \otimes V_{r}}
(W_{1}\otimes \cdots \otimes W_{r})=\C.$$
\ep

{\bf Proof.} In view of Lemma \ref{lvacounatble},
we have $\End_{V_{i}}W_{i}=\C$ for $i=1,\dots,r$.
In particular, Proposition holds for $r=1$.
For $r=2$, first it follows from 
Lemma \ref{ldensity} that $W_{1}\otimes W_{2}$ is an irreducible
$V_{1}\otimes V_{2}$-module. Then by
Lemma \ref{lvacounatble}, $\End_{V_{1}\otimes V_{2}}(W_{1}\otimes W_{2})=\C$.
Now, Proposition follows immediately from induction on $r$ and 
the assertion for $r=2$.
$\;\;\;\;\;\Box$

Now we are ready to prove our main result:

\bt{tmain}
Let $V$ be an irreducible vertex algebra of countable dimension and 
let $W$ be any $V$-module. 
Then for every positive integer $n$,
the linear map $Z_{n}^{W}$ is injective. 
Furthermore, every irreducible vertex algebra $V$ of countable dimension
is nondegenerate. 
\et

{\bf Proof.} In view of Proposition \ref{pkernal}, for any positive
integer $n$, $\ker Z_{n}^{W}$ is a 
$V^{\otimes (n+1)}\otimes \C((z_{1}))\dots ((z_{n}))$-submodule of
$V^{\otimes n}\otimes W\otimes \C((z_{1}))\dots ((z_{n}))$.
It follows from Proposition \ref{pgeneralizationfhl} that 
$V^{\otimes n}$ is an irreducible
$V^{\otimes n}$-module and that 
$\End _{V^{\otimes n}}V^{\otimes n}=\C$.
Let $A_{1}$ be the subalgebra of $\End V^{\otimes n}$
generated by all the operators $u_{m}$ for $u\in V^{\otimes n},\;
m\in \Z$. Then $A_{1}$ acts
irreducibly and $\End _{A_{1}}V^{\otimes n}=\C$.
In view of Lemma \ref{ldensity}, we have 
$$\ker Z_{n}^{W}=V^{\otimes n}\otimes U,$$ 
where $U$ is a $V\otimes \C((z_{1}))\dots ((z_{n}))$-submodule of 
$W\otimes \C((z_{1}))\dots ((z_{n}))$.
Let
$$F=w_{1}\otimes f_{1}+\cdots +w_{r}f_{r}$$
be a generic element of $U$, where $w_{1},\dots,w_{r}$ 
are linearly independent vectors in $W$ and
$f_{i}\in \C((z_{1}))\cdots ((z_{n}))$.
Then we have ${\bf 1}^{\otimes n}\otimes F\in \ker Z_{n}^{W}$, 
so that
\begin{eqnarray}
f_{1}w_{1}+\cdots + f_{r}w_{r}=
\sum_{i=1}^{r}f_{i}Y({\bf 1},z_{1})\cdots Y({\bf 1},z_{n})w_{i}
=Z_{n}^{W}({\bf 1}^{\otimes n}\otimes F)
=0.
\end{eqnarray}
Consequently, $f_{i}=0$ for $i=1,\dots,r$. This proves that
$U=0$, so $\ker Z_{n}^{W}=0$.
Furthermore, since $Z_{n}=Z_{n}^{V}E_{n}$, where $E_{n}$ is the
embedding of $V^{\otimes n}$ into $V^{\otimes (n+1)}$, $Z_{n}$ must be
injective. This proves that $V$ is nondegenerate.
$\;\;\;\;\Box$

It follows from the definition of the notion of vertex operator
algebra ([FLM], [FHL]) that any vertex operator algebra 
in the sense of [FLM] and [FHL] has  
countable dimension. We also know that a simple vertex operator
algebra is irreducible. Then we immediately have:

\bc{csimplevoa}
Let $V$ be a simple vertex operator algebra in the sense of
[FLM] and [FHL] and let $W$ be a $V$-module. 
Then for every positive integer $n$,
the linear map $Z_{n}^{W}$ is injective. 
Furthermore, every simple vertex operator algebra $V$
is nondegenerate.$\;\;\;\;\Box$
\ec

We also immediately have the following generalization of
Proposition \ref{pdmdl}:

\bc{cconstant}
Let $V$ be an irreducible vertex algebra of countable dimension and
let $W$ be any (not necessarily irreducible) $V$-module. 
Then for each positive integer $n$, the linear map
$$F_{n}^{W}:
V^{\otimes n}\otimes W\rightarrow W((z_{1}))\cdots ((z_{n}))$$ 
defined by
\begin{eqnarray}
F_{n}^{W}(v^{(1)}\otimes \cdots\otimes v^{(n)}\otimes w)=
Y(v^{(1)},z_{1})\cdots Y(v^{(n)},z_{n})w
\end{eqnarray}
is injective. $\;\;\;\;\Box$
\ec

\br{rintertwining}
{\em Results of this paper can be appropriately generalized
in terms of intertwining operators [FHL].
In fact, the corresponding result of [DL] was
formulated in terms of intertwining operators in the more general
context of generalized vertex algebras.}
\er


\end{document}